# Estimating a Polya frequency function$_2$

**Jayanta Kumar Pal**[1,*], **Michael Woodroofe**[1,†] **and Mary Meyer**[2,†]

*University of Michigan and University of Georgia*

**Abstract:** We consider the non-parametric maximum likelihood estimation in the class of Polya frequency functions of order two, viz. the densities with a concave logarithm. This is a subclass of unimodal densities and fairly rich in general. The NPMLE is shown to be the solution to a convex programming problem in the Euclidean space and an algorithm is devised similar to the iterative convex minorant algorithm by Jongbleod (1999). The estimator achieves Hellinger consistency when the true density is a $\mathrm{PFF}_2$ itself.

## 1. Introduction

The problem of estimating a unimodal density and its mode has attracted a wide interest in the literature, beginning with the work of Barlow [1], Prakasa Rao, [7], Robertson [8] and Wegman [15], [16] and continuing through [9], [2], [3], [5], and [4], who can be consulted for further references. Asymptotic properties of the maximum likelihood estimators have been developed but may be messy and suffer from some inconsistency in the region near the mode. Kernel estimation can avoid the inconsistency, but must confront the choice of a bandwidth. Here we investigate a smaller, easier version of the problem, estimating a Polya frequency function of order two [hereafter $\mathrm{PFF}_2$]. By $\mathrm{PFF}_2$, we mean a density $f$ whose logarithm is concave over the support of $f$. Equivalently, a function $f$ whose logarithm is concave in the sense of Rockafellar [10]. Such functions are automatically unimodal. Moreover, an estimated $\mathrm{PFF}_2$ supplies its own estimate of the mode. There is no need to estimate the mode seperately.

The non-parametric maximum likelihood estimator [hereafter, NPMLE] for this problem is derived in Section 2 and shown to be Hellinger consistent in Section 3. Simulations are reported in Section 4. Rufibach and Duembgen [11] and Walther [14] adopt a similar approach but obtain different results by different methods.

## 2. The NPMLE

Let $\mathcal{F}$ be the class of $\mathrm{PFF}_2$ densities and suppose that a sample $X_1, \ldots, X_n \sim f \in \mathcal{F}$ is available. The problem is to estimate $f$ non-parametrically. Letting $-\infty < x_1 < x_2 < \cdots < x_n < \infty$ denote the order statistics the log-likelihood function is

$$\ell(g) = \sum_{i=1}^{n} \log[g(x_i)]. \tag{2.1}$$

*Supported by Horace Rackham Graduate School and National Science Foundation.
†Supported by National Science Foundation.
[1]436 West Hall, 1085 South University, Department of Statistics, Ann Arbor, MI 48105, USA, e-mail: jpal@umich.edu; michaelw@umich.edu
[2]223 Statistics Building, University of Georgia, Athens, GA, USA, e-mail: mmeyer@stat.uga.edu

This is to be maximized with respect to $g \in \mathcal{F}$. Equivalently, this is to be maximized with respect to non-negative $g$ for which $\log(g)$ is concave and

$$\int_{-\infty}^{\infty} g(x)dx = 1. \tag{2.2}$$

If $g \in \mathcal{F}$, write $\theta_i = \log[g(x_i)]$. Then $\log[g(x)] \geq [(x - x_{i-1})\theta_i + (x_i - x)\theta_{i-1}]/(x_i - x_{i-1})$ for $x_{i-1} \leq x \leq x_i$ and, therefore,

$$\int_{x_{i-1}}^{x_i} g(x)dx \geq \big[\frac{e^{\theta_i} - e^{\theta_{i-1}}}{\theta_i - \theta_{i-1}}\big](x_i - x_{i-1}) \tag{2.3}$$

for $i = 2, \ldots, n$. *It follows easily that (2.1) is maximized when $g \in \mathcal{F}$ has support $[x_1, x_n]$ and $\log(g)$ is a piecewise linear function with knots at $x_1, \ldots, x_n$.* For if $g \in \mathcal{F}$, let $g^o$ be a function for which $\log(g^o)$ is piecewise linear, $g^o(x_i) = g(x_i),\ i = 1, \ldots, n$, and $g^o(x) = 0$ for $x \notin [x_1, x_n]$. Then, using (2.3), $g(x) \geq g^o(x)$ for all $x$ with equality for $x \in \{x_1, \ldots, x_n\}$, and therefore,

$$1 = \int_{-\infty}^{\infty} g(x)dx \geq \int_{x_1}^{x_n} g(x)dx \geq \int_{x_1}^{x_n} g^o(x)dx = \int_{-\infty}^{\infty} g^o(x)dx.$$

So, there is a $c \leq 1$ for which $g^* = g^o/c \in \mathcal{F}$ and $\ell(g) \leq \ell(g^*)$.

Thus, finding the NPMLE may be reformulated as a maximization problem in $I\!R^n$. Let $\mathcal{K}$ be the set of $\theta = [\theta_1, \theta_2, \ldots, \theta_n] \in I\!R^n$ for which

$$\frac{\theta_i - \theta_{i-1}}{x_i - x_{i-1}} \geq \frac{\theta_{i+1} - \theta_i}{x_{i+1} - x_i}$$

for $i = 2, \ldots, n-1$. The reformulated maximization problem is to maximize $\theta_1 + \cdots + \theta_n$ among $\theta \in \mathcal{K}$ subject to the constraint

$$\sum_{i=2}^{n} \big[\frac{e^{\theta_i} - e^{\theta_{i-1}}}{\theta_i - \theta_{i-1}}\big](x_i - x_{i-1}) = 1. \tag{2.4}$$

Introducing a Lagrange multiplier, it is necessary to maximize

$$\psi(\theta) = \sum_{i=1}^{n} \theta_i - \lambda \sum_{i=2}^{n} \big[\frac{e^{\theta_i} - e^{\theta_{i-1}}}{\theta_i - \theta_{i-1}}\big](x_i - x_{i-1}) \tag{2.5}$$

subject to (2.4), for appropriate $\lambda$, for $\theta \in \mathcal{K}$.

The partial derivatives of $\psi$ are

$$\begin{aligned}
\frac{\partial\psi(\theta)}{\partial\theta_1} &= 1 - \lambda\big[-\frac{e^{\theta_1}}{\theta_2 - \theta_1} + \frac{e^{\theta_2} - e^{\theta_1}}{(\theta_2 - \theta_1)^2}\big](x_2 - x_1), \\
\frac{\partial\psi(\theta)}{\partial\theta_2} &= 1 - \lambda\big[\frac{e^{\theta_2}}{\theta_2 - \theta_1} - \frac{e^{\theta_2} - e^{\theta_1}}{(\theta_2 - \theta_1)^2}\big](x_2 - x_1) \\
&\quad - \lambda\big[-\frac{e^{\theta_2}}{\theta_3 - \theta_2} + \frac{e^{\theta_3} - e^{\theta_2}}{(\theta_3 - \theta_2)^2}\big](x_3 - x_2), \\
&\cdots \\
\frac{\partial\psi(\theta)}{\partial\theta_j} &= 1 - \lambda\big[\frac{e^{\theta_j}}{\theta_j - \theta_{j-1}} - \frac{e^{\theta_j} - e^{\theta_{j-1}}}{(\theta_j - \theta_{j-1})^2}\big](x_j - x_{j-1}) \\
&\quad - \lambda\big[-\frac{e^{\theta_j}}{\theta_{j+1} - \theta_j} + \frac{e^{\theta_{j+1}} - e^{\theta_j}}{(\theta_{j+1} - \theta_j)^2}\big](x_{j+1} - x_j),
\end{aligned}$$

$$\cdots$$

and

$$\frac{\partial\psi(\theta)}{\partial\theta_n} = 1 - \lambda\big[\frac{e^{\theta_n}}{\theta_n - \theta_{n-1}} - \frac{e^{\theta_n} - e^{\theta_{n-1}}}{(\theta_n - \theta_{n-1})^2}\big](x_n - x_{n-1}).$$

At the maximizing $\theta$, $\nabla\psi(\theta)^t\mathbf{1} = 0$, so that

$$\begin{aligned} n = \lambda\big[-\frac{e^{\theta_1}}{\theta_2 - \theta_1} + \frac{e^{\theta_2} - e^{\theta_1}}{(\theta_2 - \theta_1)^2}\big](x_2 - x_1) \\ + \lambda\sum_{j=2}^{n-1}\Big\{\big[\frac{e^{\theta_j}}{\theta_j - \theta_{j-1}} - \frac{e^{\theta_j} - e^{\theta_{j-1}}}{(\theta_j - \theta_{j-1})^2}\big](x_j - x_{j-1}) \\ + \lambda\big[-\frac{e^{\theta_j}}{\theta_{j+1} - \theta_j} + \frac{e^{\theta_{j+1}} - e^{\theta_j}}{(\theta_{j+1} - \theta_j)^2}\big](x_{j+1} - x_j)\Big\} \\ + \lambda\big[\frac{e^{\theta_n}}{\theta_n - \theta_{n-1}} - \frac{e^{\theta_n} - e^{\theta_{n-1}}}{(\theta_n - \theta_{n-1})^2}\big](x_n - x_{n-1}) \end{aligned}$$

There is some cancelation here, and

$$n = \lambda\sum_{j=2}^{n}\frac{e^{\theta_j} - e^{\theta_{j-1}}}{\theta_j - \theta_{j-1}}(x_j - x_{j-1}).$$

So, if (2.4) is to be satisfied, then

$$\lambda = n. \tag{2.6}$$

Now let

$$\omega_1 = \theta_1$$
$$\omega_j = \frac{\theta_j - \theta_{j-1}}{x_j - x_{j-1}}$$

for $j = 2, \ldots, n$. Then,

$$\omega_2 \geq \cdots \geq \omega_n \tag{2.7}$$

and

$$\theta_j = \omega_1 + \sum_{i=2}^{j}(x_i - x_{i-1})\omega_i \tag{2.8}$$

for $j = 1, \ldots, n$, where an empty sum is to be interpreted as 0. Let

$$\phi(\omega) = \psi(\theta) \tag{2.9}$$

and $\Delta x_j = x_j - x_{j-1}$. Then

$$\sum_{j=1}^{n}\theta_j = n\omega_1 + \sum_{i=2}^{n}(n - i + 1)\Delta x_i\omega_i,$$

$$\begin{aligned} \sum_{i=2}^{n}\frac{e^{\theta_i} - e^{\theta_{i-1}}}{\theta_i - \theta_{i-1}}\Delta x_i &= \sum_{j=2}^{n}e^{\theta_{j-1}}\frac{e^{(x_j - x_{j-1})\omega_j} - 1}{\omega_j} \\ &= \sum_{j=2}^{n}e^{\theta_{j-1}}\rho(\Delta x_j\omega_j)\Delta x_j, \end{aligned}$$

where

$$\rho(t) = \frac{e^t - 1}{t}.$$

So,

$$\phi(\omega) = n\omega_1 + \sum_{i=2}^{n}(n-i+1)\Delta x_i\omega_i - n\sum_{j=2}^{n} e^{\theta_{j-1}}\rho(\Delta x_j\omega_j)\Delta x_j.$$

Using (2.8), it follows that at the maximizing $\theta$ (or $\omega$)

$$\frac{\partial\phi(\omega)}{\partial\omega_1} = n - n\sum_{j=2}^{n} e^{\theta_{j-1}}\rho(\Delta x_j\omega_j)\Delta x_j = 0.$$

So, we are led to the problem of maximizing $\phi(\omega)$, subject to (2.7) and $\partial\phi(\omega)/\partial\omega_1 = 0$. Again using (2.8), the latter condition may be written

$$e^{-\omega_1} = \sum_{j=2}^{n} e^{\sum_{i=2}^{j-1}\Delta x_i\omega_i}\rho(\Delta x_j\omega_j)\Delta x_j. \tag{2.10}$$

To solve this, we need a version of the iterative concave majorant algorithm, similar to that of Jongbloed [6]. We start with an initial value $\omega^0 = (\omega_1^0, \ldots, \omega_n^0)$ for which (2.7) and (2.10) are satisfied. One such choice is to assume $f$ is a normal density and estimate its mean and variance from the data. The corresponding $\omega^0$ can be computed using a scaled piecewise linear version of $\log f$. Let $k = 0$.

The idea behind our algorithm is to replace the concave function $\phi$ locally near $\omega^k$ by a quadratic form of the type

$$q(\tilde{\omega}, \omega^k) = \frac{1}{2}(\tilde{\omega} - \omega^k + \Gamma(\omega^k)^{-1}\nabla\phi(\omega^k))^t\Gamma(\omega^k)(\tilde{\omega} - \omega^k + \Gamma(\omega^k)^{-1}\nabla\phi(\omega^k))$$

where $\Gamma$ is a diagonal matrix with entries $\partial^2\phi(\omega)/\partial\omega_k^2$ and $\nabla\phi$ is the gradient vector. This maximization has a geometric solution given by the left hand slopes of the concave majorant of the data cloud: $(\sum_{i=1}^{l} d_i^r, \sum_{i=1}^{l}[d_i^r\omega_i^r + b_i^r])$, where

$$b_k^r = \frac{\partial}{\partial\omega_k}\phi(\omega)|_{\omega=\omega^r},$$

$$d_k^r = -\frac{\partial^2}{\partial\omega_k^2}\phi(\omega)|_{\omega=\omega^r}. \tag{2.11}$$

This can be also characterized explicitly as,

$$\omega_i^{r+1} = \min_{2\le j\le i}\max_{i\le k\le n}\frac{\sum_{h=j}^{k}[d_h^r\omega_h^r + b_h^r]}{\sum_{h=j}^{k} d_h^r}.$$

Finally, to satisfy (2.7),

$$\omega_1^{r+1} = -\log\Big[\sum_{j=2}^{n} e^{\sum_{i=2}^{j-1}\Delta x_i\omega_i^{r+1}}\rho(\Delta x_j\omega_j^{r+1})\Delta x_j\Big].$$

To implement this, we need to compute the partial derivatives $(\partial/\partial\omega_k)\phi(\omega)$ for $k = 1, \ldots, n$. Clearly,

$$\frac{\partial\phi(\omega)}{\partial\omega_1} = 0$$

Table 1
*Angular distances from the center R and line-of-sight velocities V of stars in Fornax.*

| Star | Date | $R$ arcmin | $V$ km/sec |
|---|---|---|---|
| F1-1 | 29 Nov. 1992 | 14.5 | 55.8 |
| F1-2 | 29 Nov. 1992 | 20.8 | 42.9 |
| ... | | ... | |
| F9-8025 | 15 Dec. 2002 | 42.3 | 69.2 |

and

$$\frac{\partial\phi(\omega)}{\partial\omega_k} = (n-k+1)(x_k - x_{k-1}) - n\sum_{j=k+1}^{n}\Delta x_k e^{\theta_{j-1}}\rho(\Delta x_j\omega_j)\Delta x_j - ne^{\theta_{k-1}}\rho'(\Delta x_k\omega_k)\Delta x_k^2$$

for $k = 2,\ldots,n$; and the second derivatives are

$$\frac{\partial^2\phi(\omega)}{\partial\omega_k^2} = -n\Big[\sum_{j=k+1}^{n}\Delta x_k^2 e^{\theta_{j-1}}\rho(\Delta x_j\omega_j)\Delta x_j + e^{\theta_{k-1}}\rho''(\Delta x_k\omega_k)\Delta x_k^3\Big]$$

for $k = 2,\ldots,n$. To achieve stability, we modify the algorithm using a line search method as follows. It is not certain that the new point $\tilde{\omega}$ will have a larger value of $\phi$. Therefore, we need to perform a binary search along the line segment joining $\omega^k$ and $\tilde{\omega}$ to get a point $\omega^{k+1}$ such that $\phi(\omega^{k+1}) > \phi(\omega^k)$. Finally, we stop the iteration when two consecutive iterates have very close $\phi$ values.
**Example**. Walker et. al. [13] have reported the line-of-sight velocities of 178 stars in the Fornax dwarf spheroidal galaxy. The nature of the data is reported in Table 1. The full data set can be found in [13]. Figure 1 displays the estimated density of line-of-sight velocity, assuming that the later is a Polya frequency function$_2$. The sharp peak at the mode is, unfortunately, an artifact of the method.

## 3. Consistency

Let $F$ be a distribution function with density $f$ and suppose throughout that

$$-\infty < \int_{\mathbb{R}} \log(f)dF < \infty. \tag{3.1}$$

Let $X_1, X_2, \cdots \sim^{ind} F$; and let $F_n^{\#}$ be the empirical distribution function.

$$F_n^{\#}(x) = \frac{\#\{k \le n : X_k \le x\}}{n}.$$

Further, let $h$ denote the Hellinger distance between densities,

$$h^2(g_1, g_2) = \int_{\mathbb{R}} (\sqrt{g}_1 - \sqrt{g}_2)^2 dx = 2\big[1 - \int_{\mathbb{R}} \sqrt{g_1 g_2}dx\big]. \tag{3.2}$$

The purpose of this section is to prove: *If $f$ is a PFF$_2$ for which (3.1) is satisfied and $\hat{f}_n$ is the maximum likelihood estimator, then* $\lim_{n\to\infty} h^2(f, \hat{f}_n) = 0$ *w.p.*1.

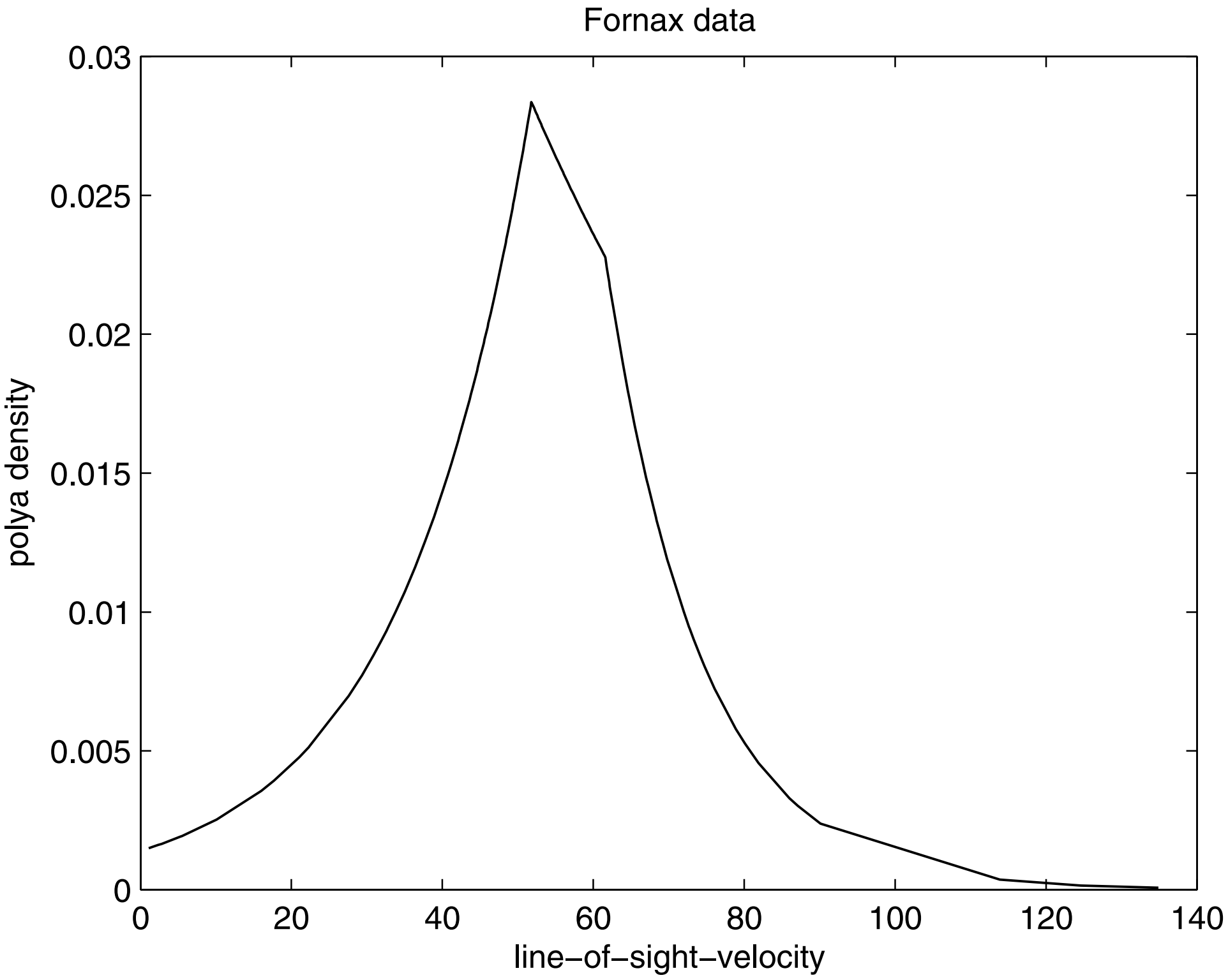

FIG 1. *Estimated density of velocities for 178 stars.*

**Lemma 1.** *If $f$ and $g$ are densities and $b > 0$, then*

$$\int_{\mathbb{R}} \log(\frac{b+g}{b+f}) dF \leq \epsilon(b) - h^2(f,g),$$

*where*

$$\epsilon(b) = 2\int_{\mathbb{R}} \sqrt{\frac{b}{b+f}} dF.$$

*Proof.* In this case,

$$\begin{aligned} \int_{\mathbb{R}} \log\left(\frac{b+g}{b+f}\right) dF &\leq 2\Big[\int_{\mathbb{R}} \sqrt{\frac{b+g}{b+f}} dF - 1\Big] \\ &\leq 2\Big[\int_{\mathbb{R}} \sqrt{\frac{b}{b+f}} dF + \int_{\mathbb{R}} \sqrt{\frac{g}{b+f}} dF - 1\Big] \\ &\leq \epsilon(b) + 2\Big[\int_{\mathbb{R}} \sqrt{gf} dx - 1\Big] \\ &= \epsilon(b) - h^2(f,g). \end{aligned}$$

□

**Lemma 2.** *Let $0 < b, c < \infty$. If $g$ is unimodal and $\sup_x g(x) \le c$, then*

$$|\int_{\mathbb{R}} \log(b+g) d(F_n^{\#} - F)| \le 2 \sup_x |F_n^{\#}(x) - F(x)| \log\big(1 + \frac{c}{b}\big). \tag{3.3}$$

*Proof.* Integrating by parts, the left side of (3.3) is

$$|\int_{\mathbb{R}} (F - F_n^{\#}) d \log(b+g)|,$$

which is at most the right side. □

Now let $\mathcal{U}$ be a class of unimodal densities; let $\ell_n$ denote the log-likelihood function, so that

$$\ell_n(g) = \sum_{i=1}^{n} \log[g(X_i)] = n \int_{\mathbb{R}} \log(g) dF_n^{\#}$$

for $g \in \mathcal{U}$; and let $\hat{f}_n$ be the MLE in $\mathcal{U}$ (assumed to exist).

**Theorem 3.1.** *If $f \in \mathcal{U}$ and $C_n = \sup_x \hat{f}_n(x)$ satisfies,*

$$\log C_n = \sup_x \log \hat{f}_n(x) = o[\frac{\sqrt{n}}{\log(n)}] \; w.p.1, \tag{3.4}$$

*then*

$$\lim_{n\to\infty} h(f, \hat{f}_n) = 0 \; w.p.1.$$

*Proof.* If $f \in \mathcal{U}$, then

$$0 \le \ell_n(\hat{f}_n) - \ell_n(f) = n\Big[\int_{\mathbb{R}} \log(\hat{f}_n) dF_n^{\#} - \int_{\mathbb{R}} \log(f) dF_n^{\#}\Big].$$

So, if $b > 0$, then

$$0 \le \int_{\mathbb{R}} \log(b + \hat{f}_n) dF_n^{\#} - \int_{\mathbb{R}} \log(f) dF_n^{\#} = I_n + II_n + III_n,$$

where

$$I_n = \int_{\mathbb{R}} \log(b + \hat{f}_n) d(F_n^{\#} - F),$$

$$II_n = \int_{\mathbb{R}} \log\big[\frac{b + \hat{f}_n}{b + f}\big] dF$$

and

$$III_n = \int_{\mathbb{R}} \log(b + f) dF - \int_{\mathbb{R}} \log(f) dF_n^{\#}.$$

With $C_n$ as in (3.4),

$$|I_n| \le 2 \sup_x |F_n^{\#}(x) - F(x)| \log\big(1 + \frac{C_n}{b}\big) \to 0 \; w.p.1$$

as $n \to \infty$, by Lemma 2 and the consistency of $F_n^{\#}$. Also,

$$II_n \le \epsilon(b) - h^2(f, \hat{f}_n),$$

by Lemma 1, and

$$\lim_{n\to\infty} III_n = \int_{I\!R} [\log(b+f) - \log(f)]dF,$$

by the Strong Law of Large Numbers. So, *w.p.*1,

$$\limsup_{n\to\infty} h^2(f, \hat{f}_n) \leq \int_{I\!R} [\log(b+f) - \log(f)]dF + \epsilon(b),$$

which approaches zero as $b \to 0$. □

**Lemma 3.** *Let $g$ be a PFF$_2$ density. If $0 < g(a) \leq g(b)$, then*

$$g(b) \leq \frac{1}{|b-a|}\Big\{1 + \log\big[\frac{g(b)}{g(a)}\big]\Big\}.$$

*Proof.* There is no loss of generality in supposing that $a < b$. Let $h = \log(g)$. Then

$$h(x) \geq h_a + (\frac{x-a}{b-a})[h_b - h_a],$$

where $h_a$ and $h_b$ were written for $h(a)$ and $h(b)$. So,

$$\begin{aligned}\int_a^b g(x)dx &\geq \int_a^b \exp\big\{h_a + (\frac{x-a}{b-a})[h_b - h_a]\big\}dx \\ &= e^{h_a}(\frac{b-a}{h_b - h_a})\big[e^{h_b - h_a} - 1\big] \\ &= (b-a)\frac{g_b - g_a}{\log(g_b/g_a)}.\end{aligned}$$

Of course, $\int_a^b g(x)dx \leq 1$, and $g(x) \geq g_a$ for $a \leq x \leq b$. So, $g_a \leq 1/(b-a)$ and

$$g_b \leq g_a + (\frac{1}{b-a})\log(\frac{g_b}{g_a}) \leq (\frac{1}{b-a})\big\{1 + \log\big[\frac{g(b)}{g(a)}\big]\big\},$$

as asserted. □

**Lemma 4.** *If $a, b, x > 0$ and $x \leq a\log(x) + b$, then $x \leq 2a\log(2a) + 2b$.*

*Proof.* If $x \leq a\log(x) + b$, then

$$x \leq a\log(2a + x) + b \leq a\big[\log(2a) + \frac{x}{2a}\big] + b = \frac{1}{2}x + a\log(2a) + b,$$

from which the lemma follows immediately. □

Now let $\hat{f}_n$ be the MLE in the class of PFF$_2$ densities. Then $\hat{f}_n$ attains its maximum at an order statistic, say $x_m$. If $m \leq n/2$, let $q = \lfloor 3n/4 \rfloor + 1$; and if $m > n/2$,let $q = \lfloor n/4 \rfloor$.

**Theorem 3.2.** *If $f$ is a PFF$_2$ density, then*

$$C := \sup_{n\geq 1} \hat{f}_n(x_m) < \infty \ w.p.1. \tag{3.5}$$

*Proof.* Let $K_n = q$ or $n - q$, accordingly as $m > n/2$ or $\leq n/2$. Then $K_n \geq n/4$, and

$$\hat{f}_n(x_m) \leq \frac{1}{|x_m - x_q|}\Big\{1 + \log\big[\frac{\hat{f}_n(x_m)}{\hat{f}_n(x_q)}\big]\Big\}, \tag{3.6}$$

by Lemma 3. If $f$ is a $\mathrm{PFF}_2$, then

$$\ell_n(f) \leq \ell_n(\hat{f}_n) \leq K_n \log[\hat{f}_n(x_q)] + (n - K_n)\log[\hat{f}_n(x_m)].$$

So,

$$\begin{aligned} \log\big[\frac{\hat{f}_n(x_m)}{\hat{f}_n(x_q)}\big] &\leq \frac{n}{K_n}\log[\hat{f}_n(x_m)] - \frac{1}{K_n}\ell_n(f) \\ &\leq 4\Big[\log[\hat{f}_n(x_m)] - \frac{1}{n}\ell_n(f)\Big]. \end{aligned} \tag{3.7}$$

Combining (3.6) and (3.7),

$$\begin{aligned} \hat{f}_n(x_m) &\leq \frac{1}{|x_m - x_q|}\Big\{1 + 4\log[\hat{f}_n(x_m)] - \frac{4}{n}\ell_n(f)\Big\} \\ &= A_n \log[\hat{f}_n(x_m)] + B_n, \end{aligned}$$

where

$$A_n = \frac{4}{|x_m - x_q|}$$

and

$$B_n = \frac{1}{|x_q - x_m|}\Big[1 - \frac{4}{n}\ell_n(f)\Big].$$

So,

$$\hat{f}_n(x_m) \leq 2A_n \log(2A_n) + 2B_n,$$

by Lemma 4. Here $\sup_n A_n < \infty$ and $\sup_n B_n < \infty$, by (3.1) and the choices of $m$ and $q$, establishing the theorem. □

From the choice of $m = m_n$, (3.5) may be written $C = \sup_n \sup_x \hat{f}_n(x)$.

**Corollary 1.** *If $\mathcal{U}$ is the class of $PFF_2$ densities and $f \in \mathcal{U}$, then* $\lim_{n\to\infty} h^2(f, \hat{f}_n) = 0$ *w.p.*1.

## 4. Simulations

To assess the speed of convergence, we conducted simulation study for different well-known members of the $\mathrm{PFF}_2$ class. The densities we sampled from are: Gaussian, Double exponential, Gamma (with shape parameter 3 and scale parameter 1), Beta (with shape parameters 3 and 2) and Weibull (with parameters 3 and 1). Table 4 gives us the summary for approximate Hellinger distances between the estimate and the true underlying density with increasing sample sizes 50, 100, 200, 500, 1000. We also graphically show how the estimators look like in some of this examples. Figure 2 shows the corresponding plots for Normal, Double Exponential and Gamma for sample sizes 50, 100, 200 respectively.

Table 2

*The Monte Carlo estimates of finite-sample Hellinger distances, for sample size $n = 50, 100, 200, 500, 1000$ and number of replications $M = 500$ for five different log-concave densities. The upper figure is the estimate and the lower is one standard deviation.*

| Sample size | Normal (0,1) | Double Exponential | Gamma (3,1) | Beta (3,2) | Weibull (3,1) |
|---|---|---|---|---|---|
| 50 | 0.1658 | 0.1548 | 0.1518 | 0.1823 | 0.0854 |
| | ± 0.0531 | ± 0.0466 | ± 0.0447 | ± 0.0582 | ± 0.0263 |
| 100 | 0.0680 | 0.1146 | 0.0934 | 0.0988 | 0.0947 |
| | ± 0.0218 | ± 0.0426 | ± 0.0305 | ± 0.0363 | ± 0.0319 |
| 200 | 0.0624 | 0.0956 | 0.0532 | 0.1166 | 0.0312 |
| | ± 0.0167 | ± 0.0252 | ± 0.0164 | ± 0.0360 | ± 0.0067 |
| 500 | 0.0290 | 0.0626 | 0.0088 | 0.0766 | 0.0347 |
| | ± 0.0104 | ± 0.0203 | ± 0.0029 | ± 0.0214 | ± 0.0115 |
| 1000 | 0.0028 | 0.0139 | 0.0019 | 0.0170 | 0.0274 |
| | ± 0.0010 | ± 0.0031 | ± 0.0007 | ± 0.0031 | ± 0.0089 |

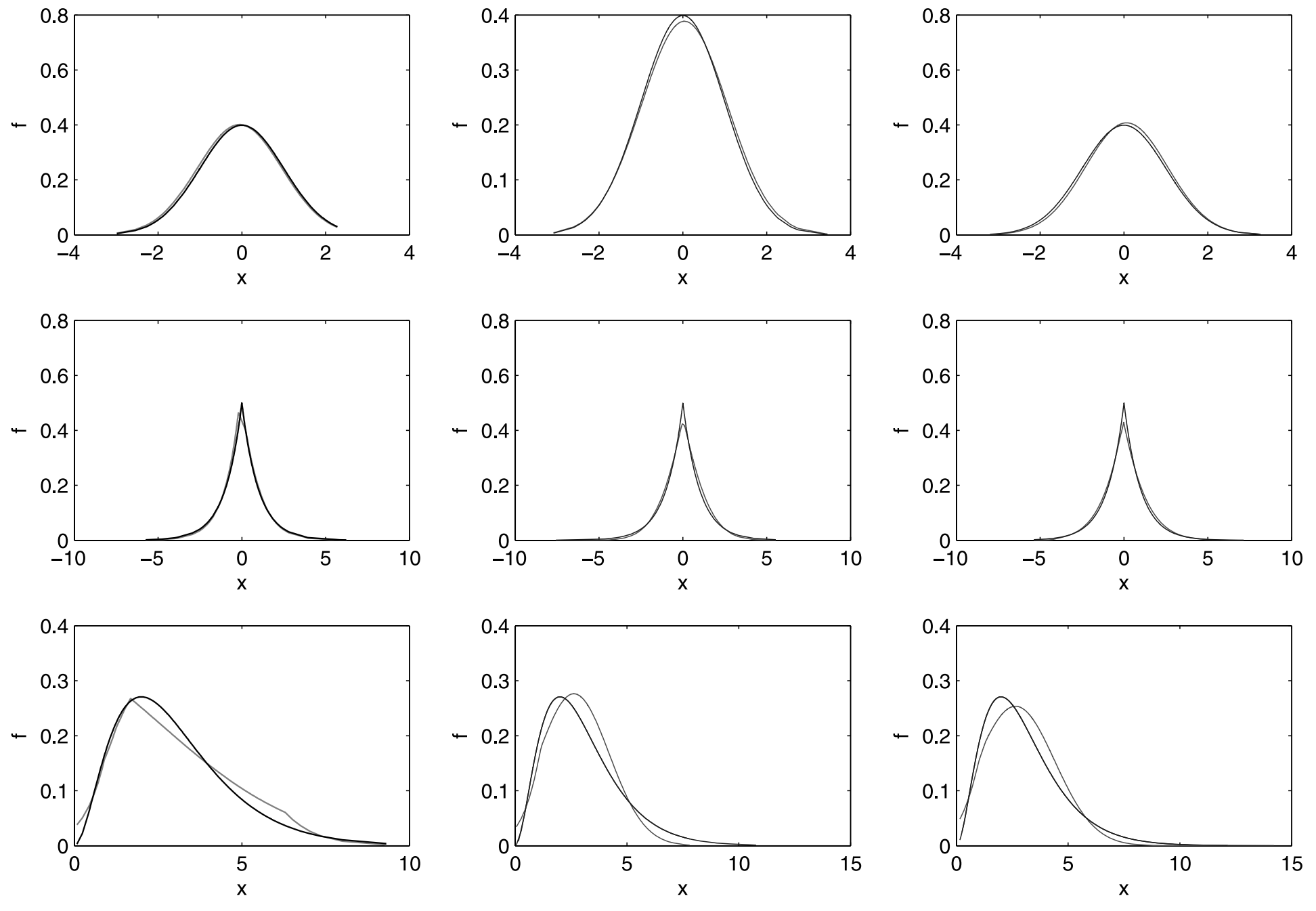

Fig 2. *The estimated log-concave density for different simulation examples. The sample sizes are 50,100 and 200 respectively for first, second and third columns. The three rows correspond to simulations from a Normal(0,1), a double-exponential and a Gamma(3,2) density. The bold one corresponds to the true density and the dotted one is the estimator.*

**Acknowledgments.** Thanks to Guenther Walther for helpful discussions. We benefited from reading Van de Geer [12] in constructing the proof of consistency.